# Semi-canonical binary matrices

**Krasimir Yordzhev**
*Faculty of Mathematics and Natural Sciences*
*South-West University, Blagoevgrad, Bulgaria*
E-mail: yordzhev@swu.bg

***Abstract***: *In this paper, we define the concepts of semi-canonical and canonical binary matrix. Strictly mathematical, we prove the correctness of these definitions. We describe and we implement an algorithm for finding all $n \times n$ semi-canonical binary matrices taking into account the number of 1 in each of them. This problem relates to the combinatorial problem of finding all pairs of disjoint $n^2 \times n^2$ S-permutation matrices. In the described algorithm, the bitwise operations are substantially used.*

## 1. INTRODUCTION

*Binary* (or *boolean*, or *(0,1)-matrix*) is called a matrix whose elements belong to the set $\mathcal{B} = \{0, 1\}$.

Let $n$ and $m$ be positive integers. With $\mathcal{B}_{n \times m}$ we will denote the set of all $n \times m$ binary matrices, while with $\mathcal{B}_n = \mathcal{B}_{n \times n}$ we will denote the set of all square $n \times n$ binary matrices.

A square binary matrix is called a *permutation matrix*, if there is just one 1 in every row and every column. Let us denote by $\mathcal{P}_n$ the group of all $n \times n$ permutation matrices, and by $\mathcal{S}_n$ the symmetric group of order $n$, i.e. the group of all one-to-one mappings of the set $[n] = \{1, 2, \ldots n\}$ in itself. In effect is the isomorphism $\mathcal{P}_n \cong \mathcal{S}_n$.

As it is well known [4,5] the multiplication of an arbitrary real or complex matrix $A$ from the left with a permutation matrix (if the multiplication is possible) leads to dislocation of the rows of the matrix $A$, while the multiplication of $A$ from the right with a permutation matrix leads to the dislocation of the columns of $A$.

Let $n$ be a positive integer and let $A \in \mathcal{B}_n$ be a $n^2 \times n^2$ binary matrix. With the help of $n-1$ horizontal lines and $n-1$ vertical lines $A$ has been separated into $n^2$ of number non-intersecting $n \times n$ square sub-matrices $A_{kl}$, $1 \le k, l \le n$, e.i.

$$A = \begin{bmatrix} A_{11} & A_{12} & \cdots & A_{1n} \\ A_{21} & A_{22} & \cdots & A_{2n} \\ \vdots & \vdots & \ddots & \vdots \\ A_{n1} & A_{n2} & \cdots & A_{nn} \end{bmatrix}.$$

A matrix $A \in \mathcal{B}_{n^2}$ is called an *S-permutation* if in each row, each column, and each sub-matrice $A_{kl}$, $1 \le k, l \le n$ of $A$ there is exactly one 1. Two S-permutation matrices $A$ and $B$ will be called *disjoint*, if there are not $i, j \in \left[ n^2 \right] = \{1, 2, \ldots, n^2\}$ such that for the elements $a_{ij} \in A$ and $b_{ij} \in B$ the condition $a_{ij} = b_{ij} = 1$ is satisfied.

The concept of S-permutation matrix was introduced by Geir Dahl [1] in relation to the popular Sudoku puzzle. Obviously a square $n^2 \times n^2$ matrix $M$ with elements of $\left[ n^2 \right] = \{1, 2, \ldots, n^2\}$ is a Sudoku matrix if and only if there are S-permutation matrices $A_1, A_2, \ldots, A_{n^2}$, each two of them are disjoint and such that $M$ can be given in the following way:

(1) $$M = 1 \cdot A_1 + 2 \cdot A_2 + \cdots + n^2 \cdot A_{n^2}.$$

In [2] Roberto Fontana offers an algorithm which randomly gets a family of $n^2 \times n^2$ mutually disjoint S-permutation matrices, where $n = 2, 3$. In $n = 3$ he ran the algorithm 1000 times and found 105 different families of nine mutually disjoint S-permutation matrices. Then using (1) he obtained $9! \cdot 105 = 38\,102\,400$ Sudoku matrices.

*Bipartite graph* is the ordered triplet $g = \langle R_g, C_g, E_g \rangle$, where $R_g$ and $C_g$ are non-empty sets such that $R_g \cap C_g = \varnothing$, the elements of which will be called *vertices*. $E_g \subseteq R_g \times C_g = \{\langle r, c \rangle \mid r \in R_g, c \in C_g\}$ - the set of *edges*. Repeated edges are not allowed in our considerations.

If $x \in \{1, 2, \ldots, n\}$, $\rho \in \mathcal{S}_n$, then the image of the element $x$ in the mapping $\rho$ we denote by $\rho(x)$. Let $g' = \langle R_{g'}, C_{g'}, E_{g'} \rangle$ and

$g'' = \langle R_{g''}, C_{g''}, E_{g''} \rangle$. We will say that the graphs $g'$ and $g''$ are *isomorphic* and we will write $g' \cong g''$, if $R_{g'} \cong R_{g''}$, $C_{g'} \cong C_{g''}$, $|R_{g'}| = |R_{g''}| = m$, $|C_{g'}| = |C_{g''}| = n$ and there exist $\rho \in \mathcal{S}_m$ and $\sigma \in \mathcal{S}_n$ such that $\langle r, c \rangle \in E_{g'} \Leftrightarrow \langle \rho(r), \sigma(c) \rangle \in E_{g''}$. In this paper we consider only bipartite graphs up to isomorphism.

Analyzing the works of G. Dahl [1] and R. Fontana [2], the question of finding a general formula for counting disjoint pairs of $n^2 \times n^2$ S-permutation matrices as a function of the integer $n$ naturally arises. This is an interesting combinatorial problem that deserves its consideration. The work [7] solves this problem. To do that, the graph theory techniques have been used. It has been shown that to count the number of disjoint pairs of $n^2 \times n^2$ S-permutation matrices, it is sufficient to obtain some numerical characteristics of the set of all bipartite graphs considered to within isomorphism of the type $g = \langle R_g, C_g, E_g \rangle$, where $V = R_g \cup C_g$ is the set of vertices, and $E_g$ is the set of edges of the graph $g$, $R_g \cap C_g = \emptyset$, $|R_g| = |C_g| = n$, $|E_g| = k$, $k = 0, 1, \ldots, n^2$.

Let $g = \langle R_g, C_g, E_g \rangle$ be a bipartite graph, where $R_g = \{r_1, r_2, \ldots, r_n\}$ and $C_g = \{c_1, c_2, \ldots, c_n\}$. Then we build the matrix $A = [a_{ij}] \in \mathcal{B}_n$, such that $a_{ij} = 1$ if and only if $\langle r_i, c_j \rangle \in E_g$. Inversely, let $A = [a_{ij}] \in \mathcal{B}_n$. We denote the $i$-th row of $A$ with $r_i$, while the $j$-th column of $A$ with $c_j$. Then we build the bipartite graph $g = \langle R_g, C_g, E_g \rangle$, where $R_g = \{r_1, r_2, \ldots, r_n\}$, $C_g = \{c_1, c_2, \ldots, c_n\}$ and there exists an edge from the vertex $r_i$ to the vertex $c_j$ if and only if $a_{ij} = 1$. It is easy to see that if $g$ and $h$ are two isomorphic graphs and $A$ and $B$ are the corresponding matrices, then $A$ is obtained from $B$ by a permutation of columns and/or rows.

Thus, the combinatorial problem to obtain and enumerate all of $n \times n$ binary matrices up to a permutation of columns or rows having exactly $k$ units naturally arises. The present work is devoted to this problem.

## 2. SEMI-CANONICAL AND CANONICAL BINARY MATRICES

**Definition 1.** Let $A \in \mathcal{B}_{n \times m}$. With $r(A)$ we will denote the ordered $n$-tuple

$$r(A) = \langle x_1, x_2, \ldots, x_n \rangle,$$

where $0 \leq x_i \leq 2^m - 1$, $i = 1, 2, \ldots n$ and $x_i$ is a natural number written in binary notation with the help of the $i$-th row of $A$.

Similarly with $c(A)$ we will denote the ordered $m$-tuple

$$c(A) = \langle y_1, y_2, \ldots, y_m \rangle,$$

where $0 \leq y_j \leq 2^n - 1$, $j = 1, 2, \ldots m$ and $y_j$ is a natural number written in binary notation with the help of the $j$-th column of $A$.

We consider the sets:

$$\begin{aligned} \mathcal{R}_{n \times m} &= \left\{ \langle x_1, x_2, \ldots, x_n \rangle \mid 0 \leq x_i \leq 2^m - 1, i = 1, 2, \ldots n \right\} \\ &= \left\{ r(A) \mid A \in \mathcal{B}_{n \times m} \right\} \end{aligned}$$

and

$$\begin{aligned} \mathcal{C}_{n \times m} &= \left\{ \langle y_1, y_2, \ldots, y_m \rangle \mid 0 \leq y_j \leq 2^n - 1, j = 1, 2, \ldots m \right\} \\ &= \left\{ c(A) \mid A \in \mathcal{B}_{n \times m} \right\} \end{aligned}$$

With "<" we will denote the lexicographic orders in $\mathcal{R}_{n \times m}$ and in $\mathcal{C}_{n \times m}$

It is easy to see that Definition 1 describes two mappings:

$$r : \mathcal{B}_{n \times m} \to \mathcal{R}_{n \times m}$$

and

$$c : \mathcal{B}_{n \times m} \to \mathcal{C}_{n \times m},$$

which are bijective and therefore

$$\mathcal{R}_{n \times m} \cong \mathcal{B}_{n \times m} \cong \mathcal{C}_{n \times m}.$$

**Definition 2.** Let $A \in \mathcal{B}_{n \times m}$,

$$r(A) = \langle x_1, x_2, \ldots, x_n \rangle,$$

$$c(A) = \langle y_1, y_2, \ldots, y_m \rangle.$$

We will call the matrix $A$ *semi-canonical*, if

$$x_1 \leq x_2 \leq \cdots \leq x_n$$

and

$$y_1 \leq y_2 \leq \cdots \leq y_m.$$

**Proposition 1.** Let $A = [a_{ij}] \in \mathcal{B}_{n \times m}$ be a semi-canonical matrix. Then there exist integers $i, j$, such that $1 \leq i \leq n$, $1 \leq j \leq m$ and

(2) $\qquad a_{11} = a_{12} = \cdots = a_{1j} = 0, \quad a_{1\,j+1} = a_{1\,j+2} = \cdots = a_{1m} = 1,$

(3) $\qquad a_{11} = a_{21} = \cdots = a_{i\,1} = 0, \quad a_{i+1\,1} = a_{i+2\,1} = \cdots = a_{n\,1} = 1.$

Proof. Let $r(A) = \langle x_1, x_2, \ldots x_n \rangle$ and $c(A) = \langle y_1, y_2, \ldots y_m \rangle$. We assume that there exist integers $p$ and $q$, such that $1 \leq p < q \leq m$, $a_{1p} = 1$ and $a_{1q} = 0$. In this case $y_p > y_q$, which contradicts the condition for semi-canonicity of the matrix $A$. We have proven (2). Similarly, we prove (3) as well. □

**Corollary 1.** Let $A = [a_{ij}] \in \mathcal{B}_{n \times m}$ be a semi-canonical matrix. Then there exist integers $s, t$, such that $0 \leq s \leq m$, $0 \leq t \leq n$, $x_1 = 2^s - 1$ and $y_1 = 2^t - 1$. □

**Definition 3.** Let $A, B \in \mathcal{B}_{n \times m}$. We will say that the matrices $A$ and $B$ are equivalent and we will write

(4) $\qquad\qquad\qquad A \sim B,$

if there exist permutation matrices $X \in \mathcal{P}_n$ and $Y \in \mathcal{P}_m$, such that

(5) $\qquad\qquad\qquad A = XBY.$

In other words $A \sim B$ if $A$ is received from $B$ after dislocation of some of the rows and the columns of $B$.

Obviously, the introduced relation is an equivalence relation.

**Definition 4.** We will call the matrix $A \in \mathcal{B}_{n \times m}$ canonical matrix, if $r(A)$ is a minimal element about the lexicographic order in the set $\{r(B) \mid B \sim A\}$.

If the matrix $A \in \mathcal{B}_{n \times m}$ is canonical and $r(A) = \langle x_1, x_2, \ldots, x_n \rangle$, then obviously

(6) $\qquad\qquad\qquad x_1 \leq x_2 \leq \cdots \leq x_n.$

From definition 4 immediately follows that in every equivalence class about the relation "$\sim$" (definition 3) there exists only one canonical matrix. Therefore, to find all bipartite graphs of type $g = \langle R_g, C_g, E_g \rangle$, where

$V = R_g \cup C_g$ is the set of vertices, and $E_g$ is the set of edges of the graph $g$, $R_g \cap C_g = \emptyset$, $|R_g| = |C_g| = n$, $|E_g| = k$, up to isomorphism, it suffices to find all canonical matrices with $k$ 1's from the set $\mathcal{B}_{n \times n}$.

With $\mathcal{T}_n \subset \mathcal{P}_n$ we denote the set of all *transpositions* in $\mathcal{P}_n$, i.e. the set of all $n \times n$ permutation matrices, which multiplying from the left an arbitrary $n \times m$ matrix swaps the places of exactly two rows, while multiplying from the right an arbitrary $k \times n$ matrix swaps the places of exactly two columns.

**Theorem 1.** *Let $A$ be an arbitrary matrix from $\mathcal{B}_{n \times m}$. Then:*

a) *If $X_1, X_2, \cdots, X_s \in \mathcal{T}_n$ are such that*

$$r(X_1 X_2 \ldots X_s A) < r(X_2 X_3 \ldots X_s A) < \cdots < r(X_s A) < r(A),$$

*then*

$$c(X_1 X_2 \ldots X_s A) < c(A).$$

b) *If $Y_1, Y_2, \cdots, Y_t \in \mathcal{T}_m$ are such that*

$$c(A Y_1 Y_2 \ldots Y_t) < c(A Y_2 Y_3 \ldots Y_t) < \cdots < c(A X_t) < r(A),$$

*then*

$$r(A Y_1 Y_2 \ldots Y_t) < r(A).$$

Proof. a) Induction by $s$. Let $s = 1$ and let $X \in \mathcal{T}_n$ be a transposition which multiplying an arbitrary matrix $A = [a_{ij}] \in \mathcal{B}_{n \times m}$ from the left swaps the places of the rows of $A$ with numbers $u$ and $v$ ($1 \le u < v \le n$), while the remaining rows stay in their places. In other words if

$$A = \begin{bmatrix} a_{11} & a_{12} & \cdots & a_{1r} & \cdots & a_{1m} \\ a_{21} & a_{22} & \cdots & a_{2r} & \cdots & a_{2m} \\ \vdots & \vdots & & \vdots & & \vdots \\ a_{u1} & a_{u2} & \cdots & a_{ur} & \cdots & a_{um} \\ \vdots & \vdots & & \vdots & & \vdots \\ a_{v1} & a_{v2} & \cdots & a_{vr} & \cdots & a_{vm} \\ \vdots & \vdots & & \vdots & & \vdots \\ a_{n1} & a_{n2} & \cdots & a_{nr} & \cdots & a_{nm} \end{bmatrix}$$

then

$$XA = \begin{bmatrix} a_{11} & a_{12} & \cdots & a_{1r} & \cdots & a_{1m} \\ a_{21} & a_{22} & \cdots & a_{2r} & \cdots & a_{2m} \\ \vdots & \vdots & & \vdots & & \vdots \\ a_{v1} & a_{v2} & \cdots & a_{vr} & \cdots & a_{vm} \\ \vdots & \vdots & & \vdots & & \vdots \\ a_{u1} & a_{u2} & \cdots & a_{ur} & \cdots & a_{um} \\ \vdots & \vdots & & \vdots & & \vdots \\ a_{n1} & a_{n2} & \cdots & a_{nr} & \cdots & a_{nm} \end{bmatrix},$$

where $a_{ij} \in \{0,1\}$, $1 \leq i \leq n$, $1 \leq j \leq m$.

Let $r(A) = \langle x_1, x_2, \ldots, x_u, \ldots, x_v, \ldots, x_n \rangle$. Then $r(XA) = \langle x_1, x_2, \ldots, x_v, \ldots, x_u, \ldots, x_n \rangle$. Since $r(XA) < r(A)$, then according to the properties of the lexicographic order $x_v < x_u$. According to Definition 1 the representation of $x_u$ and $x_v$ in binary notation with an eventual addition if necessary with unessential zeros in the beginning is respectively as follows:

$$x_u = a_{u1} a_{u2} \cdots a_{um},$$

$$x_v = a_{v1} a_{v2} \cdots a_{vm}.$$

Since $x_v < x_u$, then there exists an integer $r \in \{1, 2, \ldots, m\}$, such that $a_{uj} = a_{vj}$ when $j < r$, $a_{ur} = 1$ and $a_{vr} = 0$.

Hence if $c(A) = \langle y_1, y_2, \ldots, y_m \rangle$, $c(XA) = \langle z_1, z_2, \ldots, z_m \rangle$, then $y_j = z_j$ when $j < r$, while the representation of $y_r$ and $z_r$ in binary notation with an eventual addition if necessary with unessential zeroes in the beginning is respectively as follows:

$$y_r = a_{1r} a_{2r} \cdots a_{u-1,r} a_{ur} \cdots a_{vr} \cdots a_{nr},$$

$$z_r = a_{1r} a_{2r} \cdots a_{u-1,r} a_{vr} \cdots a_{ur} \cdots a_{nr}.$$

Since $a_{ur} = 1$, $a_{vr} = 0$, then $z_r < y_r$, whence it follows that $c(XA) < c(A)$.

We assume that for every $s$-tuple of transpositions $X_1, X_2, \ldots, X_s \in \mathcal{T}_n$ and for every matrix $A \in \mathcal{B}_{n \times m}$ from

$$r(X_1 X_2 \ldots X_s A) < r(X_2 \cdots X_s A) < \cdots < r(X_s A) < r(A)$$

it follows that

$$c(X_1 X_2 \ldots X_s A) < c(A)$$

and let $X_{s+1} \in \mathcal{T}_n$ be such that

$$r(X_1 X_2 \ldots X_s X_{s+1} A) < r(X_2 \cdots X_s X_{s+1} A) < \cdots < r(X_{s+1} A) < r(A).$$

According to the induction assumption $c(X_{s+1} A) < c(A)$.

We put

$$A_1 = X_{s+1} A.$$

According to the induction assumption from

$$r(X_1 X_2 \ldots X_s A_1) < r(X_2 \cdots X_s A_1) < \cdots < r(X_s A_1) < r(A_1)$$

it follows that

$$c(X_1 X_2 \cdots X_s X_{s+1} A) = c(X_1 X_2 \cdots X_s A_1) < c(A_1) = c(X_{s+1} A) < c(A),$$

with which we have proven a).

b) is proven similarly to a).

□

Obviously in effect is also the dual to Theorem 1 statement, in which everywhere instead of the sign "<" we put the sign ">".

**Corollary 2.** *If the matrix $A \in \mathcal{B}_{n \times m}$ is a canonical matrix, then it is a semi-canonical matrix.*

Proof. Let $A \in \mathcal{B}_{n \times m}$ be a canonical matrix and $r(A) = \langle x_1, x_2, \ldots, x_n \rangle$. Then from (6) it follows that $x_1 \leq x_2 \leq \cdots \leq x_n$. Let $c(A) = \langle y_1, y_2, \ldots, y_m \rangle$. We assume that there are $s$ and $t$ such that $s \leq t$ and $y_s > y_t$. Then we swap the columns of numbers $s$ and $t$. Thus we obtain the matrix $A' \in \mathcal{B}_{n \times m}$, $A' \neq A$. Obviously $c(A') < c(A)$. From Theorem 1 it follows that $r(A') < r(A)$, which contradicts the minimality of $r(A)$.

□

In the next example, we will see that the opposite statement of Corollary 2 is not always true.

**Example 1.** We consider the matrices:

$$A = \begin{bmatrix} 0 & 0 & 1 & 1 \\ 0 & 0 & 1 & 1 \\ 0 & 1 & 0 & 0 \\ 1 & 0 & 0 & 0 \end{bmatrix} \text{ and } B = \begin{bmatrix} 0 & 0 & 0 & 1 \\ 0 & 1 & 1 & 0 \\ 0 & 1 & 1 & 0 \\ 1 & 0 & 0 & 0 \end{bmatrix}.$$

After immediate verification, we find that $A \sim B$. Furthermore $r(A) = \langle 3,3,4,8 \rangle$, $c(A) = \langle 1,2,12,12 \rangle$, $r(B) = \langle 1,6,6,8 \rangle$, $c(B) = \langle 1,6,6,8 \rangle$. So $A$ and $B$ are two equivalent to each other semi-canonical matrices, but they are not canonical. Canonical matrix in this equivalence class is the matrix

$$C = \begin{bmatrix} 0 & 0 & 0 & 1 \\ 0 & 0 & 1 & 0 \\ 1 & 1 & 0 & 0 \\ 1 & 1 & 0 & 0 \end{bmatrix},$$

where $r(C) = \langle 1,2,12,12 \rangle$, $c(C) = \langle 3,3,4,8 \rangle$.

From example 1 immediately follows that in a given equivalence class it is possible to exist more than one semi-canonical element.

## 3. ROGRAMME CODE OF AN ALGORITHM FOR FINDING ALL SEMI-CANONICAL MATRICES

Corollary 2 is useful that it is enough to seek canonical matrices from among the semi-canonical.

In this section, we are going to suggest an algorithm (Algorithm 1) for finding the semi-canonical matrices without checking all elements of the set $\mathcal{B}_{n \times m}$, described with the help of programming language C++. In the described algorithm, bitwise operations are substantially used. In [3] and [6] we prove that the representation of the elements of $\mathcal{B}_n$ using ordered n-tuples of natural numbers and bitwise operations leads to making a fast and saving memory algorithms. Similar techniques we used in the article [8], where we describe an algorithm for solving the combinatorial problem for finding the semi-canonical matrices in the set consisting of all $n \times n$ binary matrices having exactly $k$ 1's in every row and every column. The results of this work are given in the Encyclopedia of Integer Sequences [10], respectively under the numbers A229161, A229162, A229163 and A229164. N. J. A. Sloane, who cites the work [9], presents all of them.

**Algorithm 1.** Receives all $n \times n$ semi-canonical binary matrices.

```cpp
#define n …
/*
The function check(int[], int) verifies whether obtained n-tuple represents
a semi-canonical matrix and returns the number of 1's in the matrix.
*/
int check(int x[], int n)
{
    int k=0; //  The number of 1's in the matrix. If the matrix is not
             // semi-canonical, the function returns -1.
    int yj;  // The number represents the (n-j)-th column of the matrix
    int y0=-1; // The number before yi
    for (int j=n-1; j>=0; j--)
    {
       yj=0;
       for (int i=0; i<n; i++)
       {
          if (1<<j & x[i])
          {
            yj |= 1 << (n-1-i);
            k++;
          }
       }
       if (yj<y0 ) return -1;
       y0 = yj;
    }
    // This n-tuple represents a semi-canonical matrix. We print it.
    for (int i=0; i<n; i++) cout<<x[i]<<"  ";
    cout<<"k="<<k<<'\n';
    return k;
}

int main(int argc, char *argv[])
{
int x[n]; // x[n]-ordered n-tuple of integers that represent the rows
          // of the matrix
int k[n*n+1]; //k[i] - Number of semi-canonical matrices with
              // exactly i 1's, 0<= i<=n*n
int m=n*n;
   for (int i=0; i<=m; i++) k[i]=0;
   int xmax=(1<<n)-1;
   int p,c;
   for (int s=0; s<=n; s++)
   {
```

```
      for (int i=0; i<n; i++)  x[i] = (1<<s)-1;
      c=check(x,n);
      k[s*n]++;
      p=n-1;
      while (p>0 && x[p]<xmax)
      {
        x[p]++;
        for (int i=p+1; i<n; i++)  x[i]=x[p];
        c = check(x,n);
        if (c>=0) k[c]++;
        p=n-1;
        while (x[p] == xmax)  p--;
      }
    }
  }

    for (int i=0; i<=m; i++)
        cout<<"k("<<n<<","<<i<<") = "<<k[i]<<endl;
}
```

## 4. RESULTS

Let us denote with $\kappa(n,i)$ the number of all $n \times n$ semi-canonical binary matrices with exactly $i$ 1's, where $0 \le i \le n^2$. Using Algorithm 1, we received the following integer sequences:

$$\{\kappa(2,i)\}_{i=0}^{4} = \{1,1,3,1,1\}$$

$$\{\kappa(3,i)\}_{i=0}^{9} = \{1,1,3,8,10,9,8,3,1,1\}$$

$$\{\kappa(4,i)\}_{i=0}^{16} = \{1,1,3,8,25,49,84,107,121,101,72,41,24,8,3,1,1\}$$

$$\{\kappa(5,i)\}_{i=0}^{25} = \{1,1,3,8,25,80,220,524,1057,1806,2671,3365,3680,3468,2865,$$
$$2072,1314,723,362,166,72,24,8,3,1,1\}$$

$$\{\kappa(6,i)\}_{i=0}^{36} = \{1,1,3,8,25,80,283,925,2839,7721,18590,39522,74677,125449,$$
$$188290,252954,305561,332402,326650,290171,233656,170704,113448,68677,$$
$$37996,19188,8910,3847,1588,613,299,72,24,8,3,1,1\}$$

## 5. REFERENCES

[1] Dahl, G. (2009) Permutation matrices related to sudoku, Linear Algebra and its Applications 430 (8–9), 2457–2463.